\documentclass[11pt]{article}
\usepackage{amsmath}
\usepackage{amssymb}
\usepackage{amscd}
\usepackage[all]{xy}
\usepackage{epsfig}
\newtheorem{theo}{Theorem}[section]
\newtheorem{prop}[theo]{Proposition}
\newtheorem{coro}[theo]{Corollary}
\newtheorem{lemm}[theo]{Lemma}
\newtheorem{defi}[theo]{Definition}
\newtheorem{rem}[theo]{Remark}
\newtheorem{ex}[theo]{Example}

\title{\textbf{Free Wreath Product by the Quantum 
Permutation Group}}
\date{} 
\author{\textsc{Julien Bichon}}

\makeatletter
\renewcommand{\@makefnmark}{}
\makeatother

\begin{document}

\maketitle

\hrule

\begin{abstract}
Let $A$ be a compact quantum group, let $n \in \mathbb N^*$
 and let $A_{aut}(X_n)$ be the quantum permutation group on $n$ letters.
A free wreath product
construction $A \! *_{\rm w} \! A_{aut}(X_n)$
is done. This construction provides new examples of quantum groups,
and is useful to describe the quantum automorphism group of the
$n$-times disjoint union of a finite
connected graph.
\end{abstract}

\noindent
{\small Keywords: Quantum permutation group, Wreath product, Free product,
Graph automorphism.}

\bigskip

\hrule

\section{Introduction}

We discuss a quantum analogue of the following well-known
construction. Let $n \in \mathbb N^*$ and let $G$ be a 
subgroup of the permutation group $S_n$. Let $H$ be an arbitrary
group. Then $G$ has a natural action on $H^n$ by automorphisms, and
so we may form the semi-direct product
$H^n \rtimes G$, known as the wreath product 
of $H$ by $G$ and denoted $H{\rm w}G$.
One of the most famous examples of such a construction is 
the hyperoctaedral group 
$(\mathbb Z / 2 \mathbb Z)^n \rtimes S_n = \mathbb Z / 2 \mathbb Z 
{\rm w} S_n$, the
isometry group of a hypercube in $\mathbb R^n$. 
Wreath products also occur naturally in the study of automorphism
groups of finite graphs.

In this paper the classical permutation group $S_n$ is replaced
by the quantum permutation group $A_{aut}(X_n)$. This compact quantum
group, introduced by S. Wang \cite{[Wa]}, is the universal
compact quantum group acting on $n$ points. It is
quite different
from the early examples of compact quantum groups constructed
by 
S.L. Woronowicz \cite{[WO],[WO2]}, which where related
to the concept of ``$q$-deformation'' of a Lie group. 
The representation theory of the quantum permutation group has
been worked out by 
T. Banica \cite{[Ba]}: if $n \geq 4$, the fusion
semi-ring of $A_{aut}(X_n)$ is identical to the one of $SO(3)$.
The quantum permutation group is a fascinating object from many 
perspectives.
For example we have a paradox from the noncommutative topology viewpoint.
On one hand since the $C^*$-algebra $A_{aut}(X_n)$ is generated by 
projections, the corresponding quantum space should be totally
disconnected. On the other hand since connectedness properties 
of a compact group may be read off from its fusion semi-ring
\cite{[Il],[HR]}, Banica's classification of
representations should imply that the quantum
permutation group is connected! 

In order to find natural families of quantum subgroups of the 
quantum permutation group, the quantum automorphism group
of a finite graph was introduced in \cite{[Bi]}. It turns out
that one gets interesting quantum groups when considering
the quantum automorphism group of the $n$-times disjoint union
of a finite connected graph. 
Since the classical automorphism group of such a graph
is described by a wreath product, it was natural to think 
that an analogue of the wreath product should be helpful in the quantum case.

Let $A$ be a compact quantum group.
The free wreath product of $A$ by 
$A_{aut}(X_n)$, denoted $A \! *_{\rm w}  \! A_{aut}(X_n)$, 
is a quotient of the 
free product $C^*$-algebra $A^{*n} \! * \! A_{aut}(X_n)$ (see Sections 2-3), 
and so is not a free product, but nor is a tensor product.
This construction yields new examples of quantum groups.
When $\mathcal G$ is a finite connected graph with quantum
automorphism group denoted by $A_{aut}(\mathcal G)$, we have a 
compact quantum group isomorphism
$A_{aut}(\mathcal G^{\amalg n}) \cong  
A_{aut}(\mathcal G) \! *_{\rm w} \! A_{aut}(X_n)$, analogous
to the classical isomorphism 
Aut$(\mathcal G^{\amalg n}) \cong  
{\rm Aut}(\mathcal G){\rm w} S_n$.
In this way we get a  simpler presentation for the
algebra $A_{aut}(\mathcal G^{\amalg n})$.

\medskip

Our work is organized as follows. Section 2 is devoted to the 
algebraic construction of the free wreath product.
Some new examples of Hopf algebras are constructed, for which the 
corepresentation theory is studied in special cases. In Section 3,
the free wreath product construction is done at the Woronowicz
algebra level, using the results of the previous section.
The free wreath product is used in section 4 to get a simple 
description of the quantum automorphism group of the $n$-times disjoint union
of a finite connected graph.

\medskip

We work over the field of complex numbers. We assume the reader to be familiar
with Hopf algebras, Hopf $*$-algebras, CQG algebras
and Woronowicz algebras (compact quantum groups).
The relevant definitions may be found in the book \cite{[KS]}.

\section{The Hopf algebra construction}

Let $n \in \mathbb N^*$ be a positive integer. The main character
of this paper is Wang's quantum permutation group \cite{[Wa]}.
The corresponding Hopf algebra 
$\mathcal A_{aut}(X_n)$, denoted $\mathcal A_t(n)$ for simplicity, 
is the universal (complex) algebra with
generators  $(x_{ij})_{1\leq i,j \leq n}$ and satisfying the relations:
$$x_{ij} x_{ik} = \delta_{jk} x_{ij} \quad ; \quad
x_{ji} x_{ki} = \delta_{jk} x_{ji} \quad ; \quad
\sum_{l=1}^n x_{il} = 1 = \sum_{l=1}^n x_{li} \quad 
; \quad1 \leq i,j,k \leq n.$$ 
It is immediate to check that $\mathcal A_t(n)$ is a Hopf $*$-algebra,
with structure morphisms defined by:
$$ x_{ij}^* = x_{ij} \quad ; \quad 
\Delta(x_{ij}) = \sum_{k=1}^n x_{ik} \otimes x_{kj} \quad ; \quad
\varepsilon(x_{ij}) = \delta_{ij} \quad ; \quad
S(x_{ij}) = x_{ji} \quad ;  \quad 1\leq i,j \leq n.$$
The corepresentation $x =(x_{ij})$ is unitary, and hence it follows from
\cite{[KS]}, \S11, Theorem 27, 
that $\mathcal A_t(n)$ is a CQG algebra.

\medskip

Now let $\mathcal A$ be an arbitrary algebra. We may form the free product
$\mathcal A^{*n}$, that is the $n$-times coproduct of $\mathcal A$
(in the category of unital algebras). Denote by $\nu_i : \mathcal A
\longrightarrow \mathcal A^{*n}$, $1 \leq i \leq n$, the canonical
algebras morphisms. If furthermore $\mathcal A$ is a $*$-algebra, then
$\mathcal A^{*n}$ admits a $*$-algebra structure such that the 
$\nu_i's$ are $*$-homomorphisms (there should be no confusion between
the ``$*$'' of a $*$-algebra and the ``$*$'' of a free product).
The first basic observation is that the classical action of
the permutation group on $\mathcal A^{*n}$ extends to the quantum permutation  
group.

\begin{prop}
Let $\mathcal A$ be an algebra. Then the algebra $\mathcal A^{*n}$ is,
in a natural way, a left $\mathcal A_t(n)$-comodule algebra. The coaction
$\alpha : \mathcal A^{*n} \longrightarrow \mathcal A_t(n) \otimes
\mathcal A^{*n}$ satisfies:
$$\alpha(\nu_i(a)) = \sum_{k=1}^n x_{ik} \otimes \nu_k(a) \ ,
\quad 1 \leq i \leq n \ , \quad a \in \mathcal A.$$
If furthermore $\mathcal A$ is a $*$-algebra, the coaction $\alpha$
is a $*$-homomorphism.
\end{prop}

\noindent
\textbf{Proof}. Let us first define, for i with $1 \leq i \leq n$,
a linear map $\alpha_i : \mathcal A \longrightarrow \mathcal A_t(n)
\otimes \mathcal A^{*n}$, by   
$\alpha_i(a) = \sum_{k=1}^n x_{ik} \otimes \nu_k(a)$.
Then $\alpha_i(1) = \sum_k x_{ik}\otimes 1 = 1 \otimes 1$, and 
$$\alpha_i(a) \alpha_i(b) = \sum_{k,k'} x_{ik}x_{ik'} \otimes
\nu_k(a) \nu_{k'}(b) = \sum_k x_{ik} \otimes \nu_k(ab) = \alpha_i(ab)
\ , \quad \forall a,b \in \mathcal A.$$
Thus each $\alpha_i$ is an algebra morphism, and the universal property
of the free product yields an algebra morphism
$\alpha : \mathcal A^{*n} \longrightarrow \mathcal A_t(n) \otimes
\mathcal A^{*n}$ satisfying the statement of the proposition.
It is immediate that $\alpha$ is a coaction and that
$\alpha$ is a $*$-homomorphism if $\mathcal A$ is a $*$-algebra. $\square$

\bigskip

We now assume $\mathcal A$ to be a Hopf algebra. Then the algebra
$\mathcal A^{*n}$ inherits a natural Hopf algebra structure such that
the canonical morphisms $\nu_i : \mathcal A \longrightarrow \mathcal A^{*n}$
are Hopf algebras morphisms. Using the coaction of Proposition 2.1,
we may want to form the semi-direct product 
$\mathcal A^{*n} \rtimes \mathcal A_t(n)$, see for example subsection
10.2.6 in \cite{[KS]}. However the third commutativity condition 
of Definition 8 in
\cite{[KS]} does not hold, unless $n \leq 3$ in which case $\mathcal A_t(n)$
is the function algebra on the symmetric group. This leads us to a 
notion of free wreath product.

\begin{defi}
Let $n \in \mathbb N^*$ and let $\mathcal A$ be a Hopf algebra.
The free wreath product 
of $\mathcal A$ by the quantum permutation group
$\mathcal A_t(n)$ is the quotient of the algebra $\mathcal A^{*n} *
\mathcal A_t(n)$ by the two-sided ideal generated by the elements:
$$\nu_k(a) x_{ki} - x_{ki} \nu_k(a) \ , \quad 1 \leq i,k \leq n \ , \quad
a \in \mathcal A.$$
The corresponding algebra is denoted by 
$\mathcal A \! *_{\rm w} \! \mathcal A_t(n)$.
\end{defi}   

We now state the main result of the section. We use Sweedler's notation
$\Delta(a) = a_{(1)} \otimes a_{(2)}$. The class of an element
of $\mathcal A^{*n} * \mathcal A_t(n)$ is still denoted by the same symbol
in $\mathcal A \! *_{\rm w} \! \mathcal A_t(n)$.

\begin{theo}
The free wreath product
$\mathcal A \! *_{\rm w} \! \mathcal A_t(n)$ admits a 
Hopf algebra structure.
Let $a \in \mathcal A$ and let $i,j \in \{1\ldots n\}$.
The comultiplication $\Delta$ satisfies:
$$\Delta(x_{ij}) = \sum_{k=1}^n x_{ik}\otimes x_{kj} \quad ; \
\Delta(\nu_i(a)) =  \sum_{k=1}^n \nu_i\otimes \nu_k(\Delta_{\mathcal A}(a)) 
(x_{ik} \otimes 1) =  \sum_{k=1}^n \nu_i(a_{(1)}) x_{ik} \otimes
\nu_k(a_{(2)}).$$  
The counit $\varepsilon$ satisfies $\varepsilon(x_{ij}) = \delta_{ij}$
and $\varepsilon(\nu_i(a)) = \varepsilon_{\mathcal A}(a)$.
The antipode $S$ satisfies: $$S(x_{ij}) = x_{ji} \quad ; \ 
S(\nu_i(a)) = \sum_{k=1}^n \nu_k(S_{\mathcal A}(a))x_{ki}.$$
If $\mathcal A$ is a Hopf $*$-algebra, then so is 
$\mathcal A \! *_{\rm w} \! \mathcal A_t(n)$
with $x_{ij}^* = x_{ij}$ and $\nu_i(a)^* = \nu_i(a^*)$. If $\mathcal A$
is a CQG algebra, then 
$\mathcal A \! *_{\rm w} \! \mathcal A_t(n)$ is also a
CQG algebra.
\end{theo}

\noindent
\textbf{Proof}. For simplicity we put $\mathcal H =
\mathcal A \! *_{\rm w} \! \mathcal A_t(n)$.
For $1 \leq i \leq n$, define linear maps $\delta_i : \mathcal A
\longrightarrow \mathcal H \otimes \mathcal H$ by
$\delta_i(a) = \sum_{k=1}^n \nu_i(a_{(1)})x_{ik} \otimes \nu_k(a_{(2)})$.
It is easy to show that $\delta_i$ is an algebra map (see the calculation
in the proof of Proposition 2.1). Define now an algebra morphism
$\delta_{n+1} : \mathcal A_t(n) \longrightarrow \mathcal H \otimes \mathcal H$
by $\delta_{n+1}(x_{ij}) = \sum_k x_{ik} \otimes x_{kj}$.
The universal property of the free product yields an algebra morphism
$\Delta_0 :  \mathcal A^{*n} * \mathcal A_t(n) \longrightarrow \mathcal H
\otimes \mathcal H$ such that $\Delta_0 \circ \nu_i = \delta_i$, 
$1 \leq i \leq n+1$. Then
\begin{align*}
\Delta_0(\nu_k(a)x_{ki}) & = 
\left(\sum_{l=1}^n \nu_k(a_{(1)}) x_{kl} \otimes \nu_l(a_{(2)}) \right) 
\left(\sum_{r=1}^n 
x_{kr} \otimes x_{ri}\right) \\
= & \sum_{l=1}^n \nu_k(a_{(1)}) x_{kl} \otimes \nu_l(a_{(2)}) x_{li}
=  \sum_{l=1}^n x_{kl} \nu_k(a_{(1)}) \otimes x_{li} \nu_l(a_{(2)}) \\
= & \Delta_0(x_{ki} \nu_k(a)).
\end{align*}
Hence $\Delta_0$ induces an algebra morphism $\Delta : \mathcal H
\longrightarrow \mathcal H \otimes \mathcal H$ satisfying the
required identity. One has
\begin{align*}
(\Delta \otimes & {\rm id}) \circ \Delta (\nu_i(a))
= \sum_{k=1}^n \Delta(\nu_i(a_{(1)})) \Delta(x_{ik}) \otimes \nu_k(a_{(2)}) \\
& = \sum_{k,l,r} \nu_i(a_{(1)}) x_{il}x_{ir} \otimes
\nu_l(a_{(2)}) x_{rk} \otimes \nu_k(a_{(3)}) \\
& = \sum_{k,l} \nu_i(a_{(1)}) x_{il} \otimes \nu_l(a_{(2)}) x_{lk}
\otimes \nu_k(a_{(3)}) = 
({\rm id} \otimes \Delta) \circ \Delta (\nu_i(a)).
\end{align*}
Thus $\Delta$ is coassociative since the elements 
$\nu_i(a)$ and $x_{ik}$ generate $\mathcal H$ as an algebra.
The construction of the counit is left to the reader.
Let us now define the antipode of $\mathcal H$. For $1 \leq i \leq n$
let $S_i : \mathcal A \longrightarrow  \mathcal H$ be defined by
$S_i(a) = \sum_{k=1}^n \nu_k(S_{\mathcal A}(a))x_{ki}$. One has
$S_i(1) = 1$ and 
\begin{align*}
S_i(a) S_i(b) & = \sum_{k,l} \nu_k(S_{\mathcal A}(a)) x_{ki}
\nu_l(S_{\mathcal A}(b)) x_{li} = 
\sum_{k,l} \nu_k(S_{\mathcal A}(a)) x_{ki} x_{li}
\nu_l(S_{\mathcal A}(b)) \\
& = \sum_k \nu_k(S_{\mathcal A}(a)) x_{ki}
\nu_k(S_{\mathcal A}(b)) =
\sum_k \nu_k(S_{\mathcal A}(ba)) x_{ki} = S_i(ba).
\end{align*}
Thus $S_i : \mathcal A \longrightarrow \mathcal H^{\rm op}$ is an algebra
morphism. Now let $S_{n+1} : \mathcal A_t(n) \longrightarrow 
\mathcal H^{\rm op}$  be the algebra morphism defined by
$S(x_{ij}) = x_{ji}$. The universal property of the free product yields
an algebra morphism $S_0 : \mathcal A^{*n} * \mathcal A_t(n) \longrightarrow \mathcal H^{\rm op}$ such that $S_0 \circ \nu_i = S_i$, $1 \leq i \leq n+1$.
We have
\begin{align*}
S_0(\nu_k(a)x_{ki}) & = x_{ik} \sum_{l=1}^n \nu_l(S_{\mathcal A}(a))x_{lk}
= \sum_{l=1}^n x_{ik}x_{lk} \nu_l(S_{\mathcal A}(a)) =
x_{ik}\nu_i(S_{\mathcal A}(a)) = \\
= \sum_{l= 1}^n & \nu_l(S_{\mathcal A}(a))x_{lk}x_{ik}
= S_0(\nu_k(a)) S_0(x_{ki}) = S_0(x_{ki} \nu_k(a)).
\end{align*}
Thus $S_0$ induces an algebra morphism 
$S : \mathcal H \longrightarrow \mathcal H^{\rm op}$. One has
\begin{align*}
m \circ  (S \otimes {\rm id}) & \circ \Delta (\nu_i(a)) =
\sum_{k,l} x_{ki} \nu_l(S(a_{(2)}) x_{li}\nu_k(a_{(2)})  
 = \sum_{k} x_{ki} \nu_k(S_{\mathcal A}(a_{(1)})a_{(2)}) \\
= & \varepsilon_{\mathcal A}(a)1 = \varepsilon(\nu_i(a))1,
\end{align*}
and since the elements $\nu_i(a)$ and $x_{ij}$ generate $\mathcal H$
as an algebra, it follows that  
$m \circ  (S \otimes {\rm id})  \circ \Delta = u \circ \varepsilon$
($u$ denotes the unit of $\mathcal H$). A similar computation shows
that $m \circ  ({\rm id} \otimes S) \circ \Delta = u \circ \varepsilon$: 
we conclude that 
$\mathcal H =
\mathcal A \! *_{\rm w} \! \mathcal A_t(n)$ is a Hopf algebra.

It is easy to check that if $\mathcal A$ is a Hopf $*$-algebra, then so is
$\mathcal H$ with the $*$-structure announced in the theorem.
Now let us assume that $\mathcal A$ is a CQG algebra. This means that
there exists a family $(u^{\lambda})_{\lambda \in \Lambda}$ of unitary 
corepresentations of $\mathcal A$ whose matrix coefficients generate
$\mathcal A$ as an algebra. 
Let $u^{\lambda} = (u^{\lambda}_{kl})_{1 \leq k,l \leq d_\lambda}$ be such a 
corepresentation. For $1 \leq i,j \leq n$, $1 \leq k,l \leq d_\lambda$, 
put $v^{\lambda}(i,k,j,l) = \nu_i(u^{\lambda}_{kl})x_{ij}$. Then
\begin{align*}
\Delta  (v^{\lambda} & (i,k,j,l))  =
\left(\sum_{r=1}^n \sum_{s=1}^{d_\lambda} \nu_i(u^{\lambda}_{ks}) x_{ir}
\otimes \nu_r(u_{sl}^{\lambda})\right)
\left(\sum_{t=1}^n x_{it} \otimes x_{tj} \right) = \\
& = \sum_{r=1}^n \sum_{s=1}^{d_\lambda}
\nu_i(u_{ks}^\lambda)x_{ir} \otimes \nu_r(u_{sl}^\lambda)x_{rj} =
\sum_{r,s} 
v^{\lambda}(i,k,r,s) \otimes v^{\lambda}(r,s,j,l) \ ; \\
 \varepsilon(v^{\lambda}(i, & k,j,l))  = \delta_{ij} \delta_{kl} \ ; \\
\sum_{r,s} v^{\lambda}(i,k,r,s) & v^{\lambda}(j,l,r,s)^*  = 
\sum_{r,s}\nu_i(u_{ks}^\lambda) x_{ir} x_{jr} \nu_j(u_{ls}^{\lambda^*})
= \delta_{ij} \sum_{s=1}^{d_\lambda}
\nu_i(u_{ks}^\lambda u_{ls}^{\lambda^*}) = \delta_{ij} \delta_{kl} \ ; \\
\sum_{r,s} v^{\lambda}(r,s,i,k)^* & v^{\lambda}(r,s,j,l) = 
\sum_{r,s}x_{ri} \nu_r(u_{sk}^{\lambda^*}) \nu_r(u_{sl}^\lambda) x_{rj}
= \delta_{kl} \sum_r x_{ri}x_{rj} = \delta_{kl} \delta_{ij}.
\end{align*}
These computations mean that $v^\lambda = (v^{\lambda}(i,k,j,l))$
is a unitary corepresentation of $\mathcal H$. 
Let $\mathcal B$ be the subalgebra of $\mathcal H$ by 
the coefficients of the unitary corepresentations $v^\lambda$, $\lambda 
\in \Lambda$, and $x= (x_{ij})$.
Then $\sum_{j=1}^n v^{\lambda}(i,k,j,l) = \nu_i(u_{kl}^\lambda)
\sum_j x_{ij} = \nu_i(u_{kl}^\lambda)$. Thus $\mathcal B$ contains the
elements $\nu_i(u_{kl}^\lambda)$, $1 \leq i \leq n$, $\lambda \in \Lambda$,
$1 \leq k,l \leq d_\lambda$ which generate the image of $\mathcal A^{*n}$
in $\mathcal H$ as an algebra: it follows immediately that 
$\mathcal B = \mathcal H$ and hence $\mathcal H$ is a CQG algebra.
$\square$ 

\begin{rem} It is clear that
we may also define the free wreath product by any quantum subgroup
of the quantum permutation group, that is a 
homomorphic quotient of $\mathcal A_t(n)$, and that one still has an 
obvious analogue of Theorem 2.3
\end{rem}

We now examine an example where $\mathcal A$ is a group algebra.

\begin{ex}
Let $G$ be a (discrete) group and let $n \in \mathbb N^*$.
Let $\mathcal A_n(G)$ be the universal algebra with generators
$a_{ij}(g)$, $1 \leq i,j \leq n$, $g \in G$, and submitted to the 
relations ($1 \leq i,j,k \leq n$ ;   $g,h \in G$):
$$a_{ij}(g) a_{ik}(h) = \delta_{jk} a_{ij}(gh) \quad ; \quad
a_{ji}(g) a_{ki}(h) = \delta_{jk} a_{ji}(gh) \quad ; \quad
\sum_{l=1}^n a_{il}(1) = 1 = \sum_{l=1}^n a_{li}(1).$$ 
Then $\mathcal A_n(G)$ is a CQG algebra, with:
$$ a_{ij}(g)^* = a_{ij}(g^{-1}) \ ; \ 
\Delta(a_{ij}(g)) = \sum_{k=1}^n a_{ik}(g) \otimes a_{kj}(g) \ ; \
\varepsilon(a_{ij}(g)) = \delta_{ij} \ ; \
S(a_{ij}(g)) = a_{ji}(g^{-1}).$$
Furthermore $\mathcal A_n(G)$ is isomorphic with
$\mathbb C[G] \! *_{\rm w} \! \mathcal A_t(n)$.
\end{ex}

\noindent
\textbf{Proof}. It can be checked directly that $\mathcal A_n(G)$ is 
a CQG algebra. Another way to proceed is to check that the algebras
$\mathcal A_n(G)$ and 
$\mathbb C[G] \! *_{\rm w} \! \mathcal A_t(n)$ are isomorphic,
and to transport the Hopf $*$-algebra structure of
$\mathbb C[G] \! *_{\rm w} \! \mathcal A_t(n)$ on 
$\mathcal A_n(G)$. We choose this second possibility.

If $g \in G$, the corresponding element in $\mathbb C[G]$ is still
denoted by $g$.
For $1 \leq i \leq n$, let $\phi_i : \mathbb C[G] \longrightarrow 
\mathcal A_n(G)$ be the linear map defined by
$\phi_i(g) = \sum_{k=1}^n a_{ik}(g)$ for $g \in G$. Then $\phi_i(1) =1$
and 
$$\phi_i(g) \phi_i(h) = \sum_{k,l} a_{ik}(g) a_{il}(h) = \sum_k a_{ik}(gh)
= \phi_i(gh) \ , \ \forall g,h \in G.$$
Hence $\phi_i$ is an algebra morphism.
Let $\phi_{n+1} : \mathcal A_t(n) \longrightarrow \mathcal A_n(G)$
be the algebra morphism defined by $\phi_{n+1}(x_{ij}) = a_{ij}(1)$.
The algebra morphisms $\phi_1, \ldots, \phi_{n+1}$ induce an algebra
morphism $\phi_0  :\mathbb C[G]^{*n} * \mathcal A_t(n) \longrightarrow
\mathcal A_n(G)$. One has
$$\phi_0(\nu_i(g) x_{ij}) = \sum_k a_{ik}(g) a_{ij}(1) = a_{ij}(g)
= \sum_k a_{ij}(1) a_{ik}(g) = \phi_0(x_{ij} \nu_i(g)).$$
Hence $\phi_0$ induces an algebra morphism
$\phi : \mathbb C[G] \! *_{\rm w} \! \mathcal A_t(n)
\longrightarrow \mathcal A_n(G)$. It is then easy to construct
an algebra morphism $\psi : \mathcal A_n(G) \longrightarrow
\mathbb C[G] \! *_{\rm w} \! \mathcal A_t(n)$ such that
$\psi(a_{ij}(g)) = \nu_i(g) x_{ij}$, and it is straightforward
to check that $\phi$ and $\psi$ are mutually inverse isomorphisms.
It is also easy to see that the transported Hopf $*$-algebra structure
on $\mathcal A_n(G)$ is the one announced. $\square$

\medskip

When $G$ is a finitely generated group, it is clear that
$\mathcal A_n(G)$ is a finitely generated algebra and hence is 
a CMQG algebra. The presentation can be improved if $G$ is a cyclic
group.

\noindent
$\bullet$ If $G = \mathbb Z$ is the infinite cyclic group,
then $\mathcal A_n(\mathbb Z)$ is isomorphic with
the universal $*$-algebra generated by elements $a_{ij}$, 
$1 \leq i,j \leq n$, and submitted to the relations:
\begin{align*}
a_{ij}^*a_{ij} =  & a_{ij}a_{ij}^* \  ; \
\sum_l^n a_{li}^*a_{li} = 1 = \sum_l^n a_{il}^* a_{il} \ ; \ 
a_{ij} a_{ik} = 0 = a_{ij} a_{ik}^* = a_{ij}^* a_{ik} = a_{ij}^* a_{ik}^*
\ ; \\
& a_{ji} a_{ki} = 0 = a_{ji} a_{ki}^* = a_{ji}^* a_{ki} = a_{ji}^* a_{ki}^*
\ , 1 \leq i,j,k \leq n, \ j \not = k.
\end{align*}
The antipode is given by $S(a_{ij}) = a_{ji}^*$.

\smallskip

\noindent
$\bullet$ If $G = \mathbb Z / p \mathbb Z$, then 
$\mathcal A_n(\mathbb Z / p \mathbb Z)$ is isomorphic with the universal
algebra  generated by elements $a_{ij}$, 
$1 \leq i,j \leq n$, and submitted to the relations:
$$\sum_{l=1}^n a_{il}^p = 1 = \sum_{l=1}^n a_{li}^p \ ;\ 
a_{ij} a_{ik} = 0 = a_{ji} a_{ki} \ , \ 
 1 \leq i,j,k \leq n, \ j \not = k.$$
The antipode satisfies $S(a_{ij}) = a_{ji}^{p-1}$ and we have
$a_{ij}^* = a_{ij}^{p-1}$.

\medskip
We now examine the easiest case, namely $n=2$. It is immediate 
that the elements of $\mathcal A^{*2}$ commute with those of
$\mathcal A_t(2) = C(\mathbb Z / 2 \mathbb Z)$ in 
$\mathcal A \! *_{\rm w} \! \mathcal A_t(2)$ and thus
as an algebra, we have 
$\mathcal A \! *_{\rm w} \! \mathcal A_t(2) =
\mathcal A^{*2} \! \otimes  C(\mathbb Z / 2 \mathbb Z)$. 
The Hopf algebra structure of Theorem 2.3 is just the one of a 
classical semi-direct product. If $h_1$ denotes the Haar measure
on $\mathcal A^{*2}$ (the Haar measure on a free product is described
in \cite{[W]} as the  free product of the Haar measures) and $h_2$ 
denotes the Haar measure on $C(\mathbb Z / 2 \mathbb Z)$, it is
easy to see that $h = h_1 \otimes h_2$ is the Haar measure on
$\mathcal A \! *_{\rm w} \! \mathcal A_t(2)$.

Let us describe the corepresentation theory of $\mathcal A_2(G)$
for any group $G$. First let us introduce some notation.
We consider the group free product $G * G$, with the canonical
morphisms still denoted $\nu_1,\nu_2 : G \longrightarrow G * G$.
The canonical involutive group automorphism of $G*G$ is denoted by 
$\tau$, with $\tau \circ \nu_1 = \nu_2$ and $\tau \circ \nu_2 = \nu_1$.
An element $x \in  G*G$ is the unit element if and only if
$\tau(x) = x$.

\begin{prop}
Let $G$ be a non-trivial group.

\noindent
1) To any element $x \in G*G \setminus \{1 \}$ corresponds a
two-dimensional irreducible corepresentation
$v_x$ of $\mathcal A_2(G)$. Two such corepresentations $v_x$ and $v_y$
are isomorphic if and only if $x=y$ or $x = \tau(y)$.
There also exists a non-trivial one-dimensional corepresentation $d$.

\noindent
2) Any non-trivial irreducible 
corepresentation of $\mathcal A_2(G)$ is isomorphic
to one of the corepresentations listed above.

\noindent
3) One has the following fusion rules ($x,y \in G*G \setminus \{1 \}$):
$$v_x \otimes v_y \cong v_{xy} \oplus v_{x\tau(y)} \ if
\ x \not = y^{-1} \ and \ x \not = \tau(y)^{-1} \ ; \ $$
$$v_x \otimes v_{x^{-1}} \cong \mathbb C \oplus d \oplus  
v_{x\tau(x)^{-1}} \ ; \ d \otimes d \cong \mathbb C \ ; \
v_x \otimes d \cong d \otimes v_x \cong v_x .$$
\end{prop} 

\noindent
\textbf{Proof}. We use the algebra identification
$\mathcal A_2(G) \cong \mathbb C[G*G] \otimes C(\mathbb Z / 2 \mathbb Z)$.
The basic tool of the proof is Woronowicz' character theory \cite{[WO1]},
which we freely use. Let $x \in G * G$, $x \not = 1$ and put
$$A_{11}(x) = x x_{11} \ , \ A_{12}(x) = x x_{12} \ , \
A_{21}(x) = \tau(x)x_{21}\ , \  A_{22}(x) = \tau(x)x_{22}.$$
A straightforward computation shows that
$\Delta (A_{ij}(x)) = \sum_k A_{ik}(x) \otimes A_{kj}(x)$
and  that $\varepsilon(A_{ij}(x)) = \delta_{ij}$, $1 \leq i,j \leq 2$
(recall that $x_{11} = x_{22}$ and $x_{12} = x_{21}$
in $C(\mathbb Z / 2 \mathbb Z)$). Let $v_x = (A_{ij}(x))$ be the corresponding
matrix corepresentation, with associated character
$\chi_x = (x + \tau(x))x_{11}$. 
Let $h$ be the Haar measure on $\mathcal A_2(G)$ : $h = h_1 \otimes h_2$,
see the notation and remark above. Then
$$h(\chi_x \chi_x^*) = h((x + \tau(x))(x^{-1} + \tau(x^{-1}))x_{11})
=h((2 + x\tau(x^{-1}) + \tau(x) x^{-1})x_{11}) = 2h(x_{11}) =1.$$
Hence the corepresentation $v_x$ is irreducible.
We have $\chi_{\tau(x)} = (\tau(x) + \tau^2(x))x_{11} = 
(\tau(x) +x)x_{11} = \chi_x$, hence $v_x \cong v_{\tau(x)}$
Conversely, let $y \in G*G$ be such that the corepresentations $v_x$ and
$v_y$ are isomorphic. Then $\chi_x = \chi_y$, i.e. 
$(x + \tau(x))x_{11} = (y + \tau(y))x_{11}$. This implies that
$x + \tau(x) = y + \tau(y)$ in $\mathbb C[G*G]$, i.e.
$x=y$ or $x =\tau(y)$.
Let $d = x_{11} - x_{12}$. Then $\Delta(d) = d \otimes d$, we have
a non-trivial one-dimensional corepresentation $d$ with
$d \otimes d \cong \mathbb C$ since $d^2 =1$.

Let us now prove part 3). Let $x,y \in G*G\backslash \{1\}$ with
$x \not = y^{-1}$ and $x \not = \tau(y)^{-1}$. Then  
$$h(  \chi_x \chi_y  \chi_{xy}^*) =
h((x + \tau(x)) (y + \tau(y)) (y^{-1}x^{-1} + 
\tau(y) ^{-1} \tau(x)^{-1})x_{11}) =1.$$
(The assumption $x \not = y^{-1}$ and $x \not = \tau(y)^{-1}$
has been used). Thus $v_{xy}$ appears once in the decomposition
into irreducibles of $v_x \otimes v_y$.
In the same way $h(\chi_x \chi_y \chi_{x \tau(y)}^*)=1$, and 
$v_{x \tau(y)}$ appears once in the decomposition of $v_x \otimes v_y$.
The corepresentations $v_{xy}$ and $v_{x\tau(y)}$ are not isomorphic
if $x,y \not = 1$, and so for obvious dimension reasons we must have
$$v_x \otimes v_y \cong v_{xy} \oplus v_{x\tau(y)}.$$
Now $h(\chi_x \chi_{x^{-1}}) = 
h((2 + x\tau(x)^{-1} + \tau(x)x^{-1})x_{11}) = 1$:
the trivial corepresentation $\mathbb C$ appears once in $v_x \otimes
v_{x^{-1}}$. Also 
$h(\chi_x \chi_{x^{-1}} d^*) = h(\chi_x \chi_{x^{-1}}) = 1$:
the corepresentation $d$ appears once in $v_x \otimes v_{x^{-1}}$.
Then a straightforward computation shows that
$h(\chi_x \chi_{x^{-1}} \chi_{x \tau(x)^{-1}}^*) = 1$, and hence
the corepresentation $v_{x\tau(x)^{-1}}$ 
appears once in  $v_x \otimes v_{x^{-1}}$.
So we have: 
$$v_x \otimes v_{x^{-1}} \cong \mathbb C \oplus d \oplus  
v_{x\tau(x)^{-1}}.$$
Finally we have $h(\chi_x \chi_x^* d^*) = h(\chi_x d^* \chi_x^*) =1$.
This shows that $v_x \otimes d \cong v_x \cong d \otimes v_x$.

The family of irreducible corepresentations
$\Lambda = \{ \mathbb C, d, (v_x)_{x \in G*G \backslash \{ 1 \}} \}$
contains the trivial representation, is stable under tensor product
and under conjugation, and its coefficients generate linearly
$\mathcal A_2(G)$. It follows from the orthogonality relations
\cite{[WO1]} that any irreducible corepresentation of $\mathcal A_2(G)$
is isomorphic with one of this family. $\square$

\medskip

We can see now if the fusion semi-ring of
$\mathcal A_2(G)$ is identical to the one of a classical group.

\begin{coro} Let $G$ be a non-trivial group. Then the fusion semi-ring
of $\mathcal A_2(G)$ is commutative if and only if $G \cong 
\mathbb Z / 2 \mathbb Z$. In this case the fusion semi-ring
is the same as the one of the orthogonal group $O(2)$.
\end{coro}

\noindent
\textbf{Proof}. Let us assume that the fusion semi-ring is commutative.
Then for all $x,y \in G*G \backslash \{ 1 \}$, we must have
$v_x \otimes v_{x^{-1}} \cong v_{x^{-1}} \otimes v_x$.
This implies that $v_{x \tau(x)^{-1}} \cong v_{x^{-1} \tau(x)}$,
and hence $x \tau(x)^{-1} = x^{-1} \tau(x)$ or 
$x \tau(x)^{-1} = \tau(x)^{-1} x$. Thus for any element
$x \in G*G \backslash \{ 1 \}$, we have the alternative
$x^{2} =1$ or $x \tau(x) = \tau(x) x$. Let $g \in G$ with $g \not = 1$.
We cannot have $\nu_1(g) \nu_2(g) = \nu_2(g) \nu_1(g)$, so
$\nu_1(g^2) = 1$, which implies $g^2=1$. If $h \in G$
and $h \not = 1$, then ($(\nu_1(g) \nu_2(h))^2 \not = 1$, so
$\nu_1(g) \nu_2(h) \nu_2(g) \nu_1(h) = \nu_2(g) \nu_1(h) \nu_1(g) \nu_2(h)$.
This implies that $gh = 1$, so $h = g^{-1} = g$. Hence $G$ is the cyclic 
group of order two.

Conversely, assume that $G = \langle g | g^2=1 \rangle$.
Let $l : G*G \longrightarrow \mathbb N$ be the function which assigns
to an element of $G*G$ (written in reduced form) its length.
It is immediate to check that $l(x) = l(y)$ if and only if
$x=y$ or $x = \tau(y)$.
Thus by Proposition 2.6 every irreducible 2-dimensional irreducible
corepresentation of $\mathcal A_2(G)$ can be labelled  as $v_i$ for 
some $i \in \mathbb N^*$. Put $v_0= \mathbb C$. The fusion rules of
$\mathcal A_2(G)$ now read:
$$v_i \otimes v_j \cong v_{i+j} \oplus v_{|i-j|} \ if
\ i \not = j \in \mathbb N^* \ ;$$
$$v_i \otimes v_{i} \cong v_0 \oplus d \oplus  
v_{2i} \ ; \ d \otimes d \cong  v_0 \ ; \
v_i \otimes d \cong d \otimes v_i \cong v_i, \ i \in \mathbb N^*.$$
We recognize the fusion rules of the orthogonal group $O(2)$. $\square$ 

\begin{rem} It can even be shown that the category of corepresentations
of $\mathcal A_2(\mathbb Z / 2\mathbb Z)$ and the category of
representations of $O(2)$ are monoidally equivalent. In particular 
the Hopf algebra $\mathcal A_2(\mathbb Z / 2\mathbb Z)$ is cotriangular.
\end{rem}

\section{The compact quantum group construction}

We first recall some basic notions.
As usual $\otimes$ stands for the minimal $C^*$-tensor product.
There should be no confusion with the algebraic tensor product.

A compact quantum group \cite{[WO1],[WO2']}, or Woronowicz algebra,
is a unital $C^*$-algebra $A$ together with a $C^*$-homomorphism
$\Delta : A \longrightarrow A \otimes A$ and a family of unitary matrices
$(u^\lambda)_{\lambda \in \Lambda}$, with $u_\lambda \in M_{d_\lambda}(A)$,
such that:

\noindent
(1) The subalgebra $\mathcal A$ generated by the entries $(u_{ij}^\lambda)$
of the matrices $(u^\lambda)_{\lambda \in \Lambda}$ is dense in $A$.

\noindent
(2) For $\lambda \in \Lambda$ and $i,j \in \{1, \ldots, d_\lambda \}$,
one has $\Delta(u_{ij}^\lambda) = \sum_{k=1}^{d_{\lambda}}
u_{ik}^\lambda \otimes u_{kj}^\lambda$.

\noindent
(3) For $\lambda \in \Lambda$, the transpose matrix $^t\!u^\lambda$ is
invertible.

\medskip

It follows from \cite{[WO1],[WO3]} that the dense subalgebra
$\mathcal A$ is uniquely determined, and is a CQG algebra.
A Woronowicz algebra is said to be full if $A$ is the
enveloping $C^*$-algebra of $\mathcal A$. All the Woronowicz algebras
solving universal problems (see e.g. \cite{[VDW],[Wa],[Bi]})
are full, and more generally the Woronowicz algebras constructed 
using Woronowicz' Tannaka-Krein duality \cite{[WO2]} are full.

\smallskip

Wang's quantum permutation group, denoted 
$A_{aut}(X_n)$, is the enveloping $C^*$-algebra of 
the Hopf $*$-algebra $\mathcal A_{t}(n)$ of the previous section, and is
thus a full Woronowicz algebra. The generators of 
$A_{aut}(X_n)$ are still denoted $x_{ij}$, $1 \leq i,j \leq n$.

Let $A$ be a $C^*$-algebra. We consider
the free product $C^*$-algebra (see \cite{[Av],[Voi]})
$A^{*n}$, that is the $n$-times coproduct of $A$
as a  $C^*$-algebra. We still denote by $\nu_i : \mathcal A
\longrightarrow \mathcal A^{*n}$, $1 \leq i \leq n$, the canonical
$*$-homomorphisms. Recall from \cite{[W]} that a free product
of Woronowicz algebras is in a natural way, a Woronowicz algebra.
The following definition is the $C^*$-algebra analogue of Definition 2.2.

\begin{defi}
Let $n \in \mathbb N^*$ and let $A$ be a Woronowicz algebra.
The free wreath product of $A$ by the quantum permutation group
$A_{aut}(X_n)$ is the quotient of the $C^*$-algebra $A^{*n} *
A_{aut}(X_n)$ by the two-sided ideal generated by the elements:
$$\nu_k(a) x_{ki} - x_{ki} \nu_k(a) \ , \quad 1 \leq i,k \leq n \ , \quad
a \in A.$$
The corresponding $C^*$-algebra is denoted by 
$A \! *_{\rm w} \! A_{aut}(X_n)$.
\end{defi}

\begin{theo}
The free wreath product
$A \! *_{\rm w} \! A_{aut}(X_n)$ admits a 
Woronowicz algebra structure.
Let $a \in \ A$ and let $i,j \in \{1\ldots n\}$.
The coproduct $\Delta$ satisfies:
$$\Delta(x_{ij}) = \sum_{k=1}^n x_{ik}\otimes x_{kj} \quad ; \
\Delta(\nu_i(a)) =  \sum_{k=1}^n \nu_i\otimes \nu_k(\Delta(a)) (x_{ik} \otimes
1).$$  
If $A$ is a full Woronowicz algebra, then so is 
$A \! *_{\rm w} \! A_{aut}(X_n)$.
\end{theo}

\noindent
\textbf{Proof.}
We put $H =
A \! *_{\rm w} \! A_{aut}(X_n)$.
Let $\mathcal A$ and $\mathcal A_t(n)$ be the dense
CQG algebras of $A$ and $A_{aut}(X_n)$ respectively.
For $1 \leq i \leq n$, define continuous linear maps $\delta_i : A
\longrightarrow H \otimes H$ by
$\delta_i(a) = \sum_{k=1}^n \nu_i\otimes \nu_k(\Delta(a)) (x_{ik} \otimes 1)$.
Then $\delta_{i |\mathcal A}$ is a $*$-homomorphism (see 
the proof of Theorem 2.3), and since $\mathcal A$ is dense in $A$,
It is clear that $\delta_i$ is a $*$-homomorphism.
Define now a $*$-homomorphism
$\delta_{n+1} : A_{aut}(X_n) \longrightarrow H \otimes H$
by $\delta_{n+1}(x_{ij}) = \sum_k x_{ik} \otimes x_{kj}$.
The universal property of the $C^*$-algebra
free product yields a $*$-homomorphism
$\Delta_0 :  A^{*n} * A_{aut}(X_n) \longrightarrow  H
\otimes H$ such that $\Delta_0 \circ \nu_i = \delta_i$, 
$1 \leq i \leq n+1$.
We know from the proof of Theorem 2.3
that $\Delta_0$ vanishes on the elements
$\nu_k(a) x_{ki} - x_{ki} \nu_k(a)$, $1 \leq i,k \leq n$,
$a \in \mathcal A$, and since $\mathcal A$ is dense in $A$ and
$\Delta_0$ is continuous, we get a $*$-homomorphism
$\Delta : H \longrightarrow  H \otimes H$: this is the coproduct
announced in the theorem.
Consider the canonical $*$-algebra map 
$\iota : \mathcal H = \mathcal A \! *_{\rm w} \! \mathcal A_t(n) 
\longrightarrow A \! *_{\rm w} \! A_{aut}(X_n) = H$.
It is clear that $\iota(\mathcal H)$ is dense in $H$ and that $\iota$
commutes with the respective coproducts. Furthermore
$\mathcal H$ is a CQG algebra by theorem 2.3, and hence it immediate
that $H =A \! *_{\rm w} \! A_{aut}(X_n)$ is a Woronowicz algebra.

Let us now assume that $A$ is a full Woronowicz algebra.
Let us show that if $B$ is a $C^*$-algebra and 
$\pi : \mathcal H = \mathcal A \! *_{\rm w} \! \mathcal A_t(n)
\longrightarrow B$ is $*$-algebra morphism, then there exists a 
(unique) $*$-homomorphism 
$\tilde{\pi} : H = A \! *_{\rm w} \! \mathcal A_{aut}(X_n)
\longrightarrow B$ such that $\tilde{\pi} \circ \iota = \pi$. This 
will prove that $H$ is the enveloping $C^*$-algebra of 
$\iota(\mathcal H)$ and thus that 
$H = A \! *_{\rm w} \! A_{aut}(X_n)$ 
is a full Woronowicz algebra.    
For $1 \leq i \leq n$, let 
$\pi_i : \mathcal A \longrightarrow B$ be the $*$-algebra morphism
defined by the composition
\begin{equation*}
\begin{CD}
\pi_i : \mathcal A @> \nu_i >> \mathcal A^{*n} @>>>
\mathcal A \! *_{\rm w} \! \mathcal A_t(n)
@> \pi >> B.
\end{CD}
\end{equation*}
Since $A$ is full, there exists a unique $*$-homomorphism
$\hat{\pi}_i : A \longrightarrow B$ such that 
$\hat{\pi}_{i|\mathcal A} = \pi_i$.
Let $\hat{\pi}_{n+1}$ be the $*$-homomorphism
$A_{aut}(X_n) \longrightarrow B$ such that 
$\hat{\pi}_{n+1|\mathcal A_t(n)} = \pi_{|\mathcal A_t(n)}$.
The universal property of the free product yields a $*$-homomorphism
$\hat{\pi} : A^{*n} \! * \! A_{aut}(X_n) \longrightarrow B$.
It then immediate from the density of
$\mathcal A^{*n} \! * \! \mathcal A_t(n)$ in
$A^{*n} \! * \! \mathcal A_{aut}(X_n)$ that $\hat{\pi}$ induces
a $*$-homomorphism $\tilde{\pi} :
A \! *_{\rm w} \! A_{aut}(X_n) \longrightarrow B$
such that 
$\tilde{\pi} \circ \iota = \pi$. $\square$

\medskip

The examples described in Section 2 may be adapted in an obvious
manner to the Woronowicz algebra setting. Let $G$ be discrete group.
We define $A_n(G)$ to be the Woronowicz algebra 
$C^*(G) \! *_{\rm w} \! A_{aut}(X_n)$, where
$C^*(G)$ denotes the enveloping $C^*$-algebra of $G$.
It is clear that $A_n(G)$ is the enveloping $C^*$-algebra of 
the CQG algebra $\mathcal A_n(G)$ of Section 2.

\section{Application to the quantum automorphism group
of a  \\ finite graph}

We use the results of the previous sections to describe the 
quantum automorphism group of the $n$-times disjoint union of a finite
connected graph. First let us recall some definitions from \cite{[Bi]}.

In this paper a finite graph $\mathcal G = (V,E)$ consists of two
finite sets $V$ (set of vertices) and $E$ (set of edges) such that
$E \subset V \times V$. The source and target maps
$s,t : E \longrightarrow V$ are the restrictions of the first and
second projections respectively. So in our conventions, we
do not allow a graph to have multiple edges. The quantum automorphism
group of a graph $\mathcal G$ was defined in \cite{[Bi]}. It is
still possible to describe the quantum automorphism group of a finite graph
with possible multiple edges, using e.g. the general construction of 
\cite{[Bi2]}. However some motivation for the construction done in \cite{[Bi]}
was to describe non-trivial quantum subgroups of 
Wang's quantum permutation group, so it was quite natural to restrict to 
graphs for which an automorphism is uniquely determined by its action on
the vertices. 

\smallskip
 
Let $\mathcal G = (V,E)$ be a finite graph with set of vertices
$V = \{1, \ldots, m \}$. The quantum automorphism of $\mathcal G$ \cite{[Bi]},
denoted $A_{aut}(\mathcal G)$,
is  defined to be the universal $C^*$-algebra with generators 
$(X_{ij})_{1 \leq i,j \leq m}$ and relations:
$$ X_{ij}^* = X_{ij} \ ; \
X_{ij} X_{ik} = \delta_{jk} X_{ij} \  ; \ 
X_{ji} X_{ki} = \delta_{jk} X_{ji} \  ; \ 
\sum_{l=1}^m X_{il} = 1 = \sum_{l=1}^m X_{li} \  
, \ 1 \leq i,j,k \leq m \leqno(4.1)$$
$$ 
X_{s(\gamma)i} X_{t(\gamma)k} = X_{t(\gamma)k} X_{s(\gamma)i} = 0 
\leqno(4.2)$$
$$
\quad \quad \quad \quad \quad \quad
X_{i s(\gamma)} X_{k t(\gamma)} = X_{k t(\gamma)} X_{i s(\gamma)}
=0 \quad, \quad \gamma \in E, \ (i,k) \not\in E$$
$$ 
X_{s(\gamma)s(\gamma')} X_{t(\gamma)t(\gamma')} =
X_{t(\gamma)t(\gamma')} X_{s(\gamma)s(\gamma')}
\quad , \ \gamma,\gamma' \in E \leqno(4.3)$$
$$ 
\sum_{\gamma' \in E} X_{s(\gamma')s(\gamma)} X_{t(\gamma')t(\gamma)}
= 1 = \sum_{\gamma'  \in E} X_{s(\gamma)s(\gamma')} X_{t(\gamma)t(\gamma')}
\quad , \quad \gamma \in E \leqno(4.4) $$
Then $A_{aut}(\mathcal G)$ is a full Woronowicz algebra, 
with coproduct defined by $\Delta(X_{ij})= \sum_k X_{ik} \otimes X_{kj}$
It is shown in \cite{[Bi]} that $A_{aut}(\mathcal G)$ is the universal
compact quantum group acting on the graph $\mathcal G$.
The non-trivial example considered there
is the quantum automorphism group of the following graph: \ \ \ \ \ \ \
\xymatrix{*=0{\bullet} \ar@/^/[rr] && *={\bullet}
\ar@/^/[ll] &
*=0{\bullet} \ar@/^/[rr] && *={\bullet}
\ar@/^/[ll]}

\bigskip

This quantum group will easily be described using Theorem 4.2
and the results of the previous sections:
it is isomorphic with $A_2(\mathbb Z / 2 \mathbb Z)$.

\medskip

Let us give a more convenient presentation for 
$A_{aut}(\mathcal G)$.

\begin{prop}
Let $\mathcal G = (V,E)$ be a finite graph with set of vertices
$V = \{1, \ldots, m \}$.
Then $A_{aut}(\mathcal G)$
is isomorphic with the universal $C^*$-algebra with generators
$(X_{ij})_{1 \leq i,j \leq m}$ and relations (4.1)-(4.3) and
$$ 
\sum_{k, (k,j) \in E} X_{ik}
= \sum_{k, (i,k) \in E} X_{kj} \quad ;  \quad 
\sum_{k, (k,j) \in E} X_{ki}
= \sum_{k, (i,k) \in E} X_{jk} \quad, \quad 
1 \leq i,j \leq m. 
\leqno(4.4)' $$
\end{prop}

\noindent
\textbf{Proof}. Let $A$ be an algebra with elements 
$(X_{ij})_{1 \leq i,j \leq m}$ satisfying relations (4.1).
Assume that relations (4.4) hold.
Let $i,j \in \{1, \ldots , m \}$. Then
\begin{align*}
\sum_{\gamma, t(\gamma) = j} & X_{is(\gamma)} = 
\sum_{\gamma, t(\gamma) = j} X_{is(\gamma)}
(\sum_{\gamma' \in E} X_{s(\gamma')s(\gamma)} X_{t(\gamma')t(\gamma)})=
\sum_{\gamma, t(\gamma) = j} 
\sum_{\gamma', s(\gamma') = i} X_{s(\gamma')s(\gamma)} X_{t(\gamma')t(\gamma)}
\\
= & \sum_{\gamma', s(\gamma')=i}
(\sum_{\gamma \in E} X_{s(\gamma')s(\gamma)} X_{t(\gamma')t(\gamma)})
X_{t(\gamma')j} = \sum_{\gamma, s(\gamma)=i} X_{t(\gamma)j}.
\end{align*}
In the same way one shows that
$\sum_{\gamma, t(\gamma)=j} X_{s(\gamma)i} = \sum_{\gamma,s(\gamma)=i}
X_{jt(\gamma)}$. Hence relations (4.4)' hold.

Conversely, assume that relations (4.4)' are fulfilled
and let $\gamma \in E$. Then
\begin{align*}
\sum_{\gamma'\in E} & X_{s(\gamma')s(\gamma)} X_{t(\gamma')t(\gamma)} =
\sum_{i=1}^m \sum_{k, (i,k) \in E}X_{is(\gamma)} X_{kt(\gamma)} = \\ 
=  & \sum_{i=1}^mX_{is(\gamma)} \sum_{k, (k,t(\gamma)) \in E} 
X_{ik} = 
\sum_{i=1}^mX_{is(\gamma)} = 1,
\end {align*}
and
\begin{align*}
\sum_{\gamma'\in E} & X_{s(\gamma)s(\gamma')} X_{t(\gamma)t(\gamma')} =
\sum_{i=1}^m \sum_{k, (i,k) \in E}X_{s(\gamma)i} X_{t(\gamma)k} = \\ 
=  & \sum_{i=1}^mX_{s(\gamma)i} \sum_{k, (k,t(\gamma)) \in E} 
X_{ki} = 
\sum_{i=1}^mX_{s(\gamma)i} = 1.
\end {align*}
Hence relations (4.4) hold. $\square$

\medskip

We now begin our study of the quantum automorphism group of 
the $n$-times disjoint union of a finite connected graph. Let us recall
some basic definitions.

A graph $\mathcal G = (V,E)$ is said to be connected if 
$\forall (i,j) \in V \times V$, there exists a finite sequence
of edges $\gamma_1, \ldots, \gamma_r$ such that:

- $s(\gamma_1) = i$ or $t(\gamma_1) = i$,

- $s(\gamma_r) = j$ or $t(\gamma_r) =j$,

- For $k \in \{1, \ldots , r-1\}$, $s(\gamma_k) = t(\gamma_{k+1})$
or $s(\gamma_k) = s(\gamma_{k+1})$ or
$t(\gamma_k) = t(\gamma_{k+1})$ or $t(\gamma_k) = s(\gamma_{k+1})$.

\smallskip
 Let $\mathcal G_1 = (V_1,E_1)$ and let $\mathcal G_2 = (V_2,E_2)$
be two finite graphs. Their
disjoint union is defined in the obvious way:
$\mathcal G_1 \amalg \mathcal G_2 = (V_1 \amalg V_2,E_1 \amalg E_2)$.

\medskip

Let $\mathcal G = (V,E)$ be a finite connected graph and let
$n \in \mathbb N^*$. It is well-known that the groups
Aut$(\mathcal G^{\amalg n})$ and 
Aut$(\mathcal G)^n \! \rtimes \! S_n = {\rm Aut}(\mathcal G){\rm w}S_n$
are isomorphic, and it is natural to expect that such a kind of isomorphism
still holds at the quantum automorphism group level. This is exactly the
motivation for our free wreath product construction:

\begin{theo}
Let $\mathcal G$ be a finite connected graph
and let $n \in \mathbb N^*$. Then we have a Woronowicz
algebras isomorphism:
$$A_{aut}(\mathcal G^{\amalg n}) \cong  
A_{aut}(\mathcal G) \! *_{\rm w} \! A_{aut}(X_n).$$
\end{theo} 

The proof of Theorem 4.2 will be completed at the end of the section.
Let us first note an immediate consequence:

\begin{coro}
Let $\mathcal G$ be a finite connected graph
with non-trivial automorphism group and let $n \in \mathbb N^*$
with $n \geq 2$. Then $A_{aut}(\mathcal G^{\amalg n})$ is an 
infinite-dimensional noncommutative and noncocommmutative
Woronowicz algebra. $\square$
\end{coro}

As a consequence, we also have a description of the quantum automorphism group
of a  disjoint union of polygonal graphs.
Let $m \in \mathbb N^*$ and let $\mathcal P_m = (V_m,E_m)$ with
$V_m = \{1, \ldots, m \}$ and $E_m = \{(1,2), \ldots, (m-1,m), (m,1) \}$
be the polygonal graph with $m$ vertices. It is easy to check that 
$A_{aut}(\mathcal P_m) \cong C^*(\mathbb Z/ m \mathbb Z)$. Combining the
results of the previous sections and Theorem 4.2, we
have $A_{aut}(\mathcal P_m^{\amalg n}) \cong A_n(\mathbb Z/ m \mathbb Z)$.
More generally, if $\mathcal G$ is a connected graph with abelian
automorphism group $G$, we have 
$A_{aut}(\mathcal G^{\amalg n}) \cong A_n(G)$. 

\medskip

We now prove some preliminary results for the proof of 
Theorem 4.2. For this we consider a family of finite connected graphs
$\mathcal G_1, \ldots , \mathcal G_n$. For $i \in \{1 , \ldots , n \}$,
the set of vertices (resp. of edges) of $\mathcal G_i$
is denoted by $V_i$ (resp. $E_i$). The next lemmas
are basic computations in $A_{aut}(\amalg_{i=1}^n \mathcal G_i)$.

Let $k,l \in \{1 , \ldots , n \}$ and let $i \in V_l$. We define the following 
element of $A_{aut}(\amalg_{i=1}^n \mathcal G_i)$:
$$P_i^{kl} = \sum_{u \in V_k} X_{ui}.$$

\begin{lemm}
For $k,l \in \{1 , \ldots , n \}$ and $i,j \in V_l$,
one has $P_i^{kl} = P_j^{kl}$. Hence we put
$P^{kl} = P_i^{kl}$, $\forall i \in V_l$.
For $k,k',l \in \{1 , \ldots , n \}$, we have
$$P^{kl} P^{k'l} = \delta_{kk'} P_{kl} \ ; \
\sum_{k=1}^n P^{kl} = 1 \ 
; \
(P^{kl})^* = P^{kl} \ ; \
\Delta(P^{kl}) = \sum_{k'=1}^n P^{kk'} \otimes P^{k'l} \ ,$$ 
and
$$\varepsilon(P^{kl}) = \delta_{kl} \ ; \ 
P^{lk} P^{lk'} = \delta_{kk'} P^{lk} \quad ; \quad
\sum_{k=1}^n P^{lk} = 1.$$
\end{lemm} 

\noindent
\textbf{Proof}. Let $i \in V_l$. Then by relations (4.1) we have
$$P_i^{kl} P_i^{kl} = \sum_{u,v \in V_k} X_{ui}X_{vi} = P_i^{kl}.
\leqno(4.5)$$
Let $(i,j) \in E_l$, then by relations (4.2)-(4.3), we have
$$P_i^{kl} P_j^{kl} = \sum_{\gamma \in E_k} 
X_{s(\gamma)i}X_{t(\gamma)j} = P_j^{kl} P_i^{kl},
\leqno(4.6)$$
and thus by relations (4.4) we have
$$\sum_{k=1}^n P_i^{kl} P_j^{kl} = 1 = \sum_{k=1}^n P_j^{kl} P_i^{kl}.
\leqno(4.7)$$
Let $(i,j) \in E_l$ and $k,k',l \in \{1, \ldots , n \}$
with $k \not = k'$. Then by relations 
(4.2) we have 
$$P_i^{kl} P_j^{k'l} = \sum_{u \in V_k} \sum_{v \in V_{k'}}
X_{ui} X_{vj} = 0 = P_j^{k'l} P_i^{kl}. \leqno(4.8)$$
Combining relations (4.7), (4.8) and (4.5), we see that
$$P_j^{kl} = P_i^{kl} P_j^{kl} = P_j^{kl} P_i^{kl} = P_i^{kl}
\ , \quad k,l \in \{1, \ldots , n \} \ , \quad (i,j) \in E_l.$$ 
The graph $\mathcal G_l$ is connected, and hence we have
$P_i^{kl} = P_j^{kl}$, $\forall i,j \in V_l$. This proves our first
claim. The first two relations of the lemma are
(4.5,4.8) and (4.7) respectively.
Let $k,l \in \{1,\ldots ,n \}$ and let $i \in V_l$. Then
$$\Delta(P^{kl}) = \Delta(P_i^{kl}) = 
\sum_{u \in V_k} \sum_{p = 1}^n \sum_{v \in V_p}
X_{uv} \otimes X_{vi} = 
\sum_{p=1}^n \sum_{v \in V_p} P_v^{kp} \otimes X_{vi} 
= \sum_{p=1}^n P^{kp} \otimes P_i^{pl}.$$
It is obvious that $(P^{kl})^* = P^{kl}$ and 
$\varepsilon(P^{kl}) = \delta_{kl}.$
The last two relations follow from the previous ones and Wang's Theorem 3.1
in \cite{[Wa]}. $\square$

\medskip

We use Lemma 4.4 to prove the following result (we retain the previous 
notations):

\begin{lemm}
Let $k,l \in \{1,\ldots ,n \}$ and let $i \in V_k$, $i' \not \in V_k$,
$j,j' \in V_l$. Then we have
$$X_{ij} X_{i'j'} = X_{i'j'} X_{ij} = 0 = 
X_{ji} X_{j'i'} = X_{j'i'} X_{ji}.$$
\end{lemm}

\noindent
\textbf{Proof}. Let $i \in V_k$ and let $j \in V_l$. Then we have
$$X_{ij}P^{kl} = X_{ij} P_j^{kl} =
\sum_{u \in V_k} X_{ij} X_{uj} = X_{ij} = P^{kl}X_{ij}.$$  
If $i' \not \in V_k$, then
$$X_{i'j}P^{kl} = X_{i'j} P_j^{kl} =
\sum_{u \in V_k} X_{i'j} X_{uj} = 0 = P^{kl}X_{i'j}.$$
Hence  if $j,j' \in V_l$, $i \in V_k$ and $i' \not \in V_k$, one has
$$X_{ij}X_{i'j'} = X_{ij}P^{kl}X_{i'j'} = 0 = X_{i'j'} X_{ij}.$$
The second identity is obtained using the antipode on the dense CQG algebra.
 $\square$  

\subsection*{Proof of Theorem 4.2}

Let $\mathcal G$ be a finite connected graph with set of vertices
$V= \{1, \ldots , m \}$. For $k \in \{1, \ldots , n \}$, we denote by
$\mathcal G_k$ the $k$-th copy of $\mathcal G$ in $\mathcal G^{\amalg n}$,
with set of vertices $V_k = \{k(1), \ldots ,k(p) \}$. The generators of
$A_{aut}(\mathcal G)$ are denoted by $(u_{ij})_{1 \leq i,j \leq m}$, and the
generators of  $A_{aut}(X_n)$ are still denoted by $(x_{kl})_{1 \leq k,l
\leq n}$. Let $\mathcal B$ be the free algebra with generators 
$X_{k(i)l(j)}$, $1 \leq i,j \leq m$, $1 \leq k,l \leq n$, and define
an algebra morphism 
$$\Phi_0 : \mathcal B \longrightarrow  
A_{aut}(\mathcal G) \! *_{\rm w} \! A_{aut}(X_n) \ , \quad
X_{k(i)l(j)} \longmapsto \nu_k(u_{ij})x_{kl}.$$
We will show that $\Phi_0$ induces a Woronowicz algebra morphism
$A_{aut}(\mathcal G^{\amalg n}) \longrightarrow
A_{aut}(\mathcal G) \! *_{\rm w} \! A_{aut}(X_n)$:
these are straightforward computations. First let us note that if $\mathcal B$
is endowed with the $*$-algebra structure such that 
$X_{k(i)l(j)}^* = X_{k(i)l(j)}$, then $\Phi_0$ is easily seen to be
a $*$-algebra map. Let $k,l,p \in \{1, \ldots, n \}$ and 
$i,j,r \in \{1, \ldots , m \}$. Then
\begin{align*}
& \bullet \ 
\Phi_0(X_{k(i)l(j)} X_{k(i)p(r)}) = \nu_k(u_{ij})x_{kl}
\nu_k(u_{ir})x_{kp} = \nu_k(u_{ij}u_{ir})x_{kl}x_{kp}
= \delta_{l(j)p(r)} \Phi_0(X_{k(i)l(j)}). \\
& \bullet \
\Phi_0(X_ {l(j)k(i)} X_{p(r)k(i)}) = 
\nu_l(u_{ji})x_{lk} \nu_p(u_{ri})x_{pk} =
\nu_l(u_{ji})x_{lk}x_{pk} \nu_p(u_{ri}) = 
\delta_{l(j)p(r)} \Phi_0(X_{l(j)k(i)}). \\
& \bullet \
\Phi_0(\sum_{k=1}^n \sum_{i=1}^m X_{k(i)l(j)}) = 
\sum_{k=1}^n \sum_{i=1}^m \nu_k(u_{ij})x_{kl} = 
\sum_{k=1}^n x_{kl} = 1 = \Phi_0(1). \\
& \bullet \ 
\Phi_0(\sum_{l=1}^n \sum_{j=1}^m X_{k(i)l(j)}) =
\sum_{l=1}^n \sum_{j=1}^m \nu_k(u_{ij})x_{kl} = \sum_{l=1}^n x_{kl} = 1
= \Phi_0(1). 
\end{align*}
Let $(i,j) \in E = E(\mathcal G)$, let $k,l',l'' \in \{1, \ldots, n \}$ and 
$i',i'' \in \{1, \ldots , m \}$ with $l' \not = l''$ or 
$(i',i'') \not \in E$. Then:
\begin{align*}
& \bullet \
\Phi_0(X_{k(i)l'(i')} X_{k(j)l''(i'')}) = 
\nu_k(u_{ii'})x_{kl'} \nu_k(u_{ji''})x_{kl''} 
=  \nu_k(u_{ii'}u_{ji''})x_{kl'}x_{kl''} = 0 = \\ 
& = \Phi_0(X_{k(j)l''(i'')} X_{k(i)l'(i')}). \\
& \bullet \
\Phi_0(X_{l'(i')k(i)} X_{l''(i'')k(j)}) =
\nu_{l'}(u_{i'i})x_{l'k} \nu_{l''k}(u_{i''j})x_{l''k} = 
\delta_{l'l''} \nu_{l'}(u_{i'i}u_{i''j})x_{l'k} = 0 = \\
& = \Phi_0(X_{l''(i'')k(j)} X_{l'(i')k(i)}).
\end{align*}
Let $(i,j),(i',j') \in E$ and let $k,l \in \{1, \ldots ,n \}$. Then
\begin{align*}
& \bullet \
\Phi_0(X_{k(i)l(i')} X_{k(j)l(j')})
= \nu_k(u_{ii'})x_{kl} \nu_k(u_{jj'})x_{kl} = 
\nu_k(u_{ii'}u_{jj'})x_{kl} = \\
& = \nu_k(u_{jj'})x_{kl} \nu_k(u_{ii'})x_{kl}=  
\Phi_0(X_{k(j)l(j')} X_{k(i)l(i')})
\end{align*}
Let $k,l \in \{1, \ldots , n \}$ and $i,j \in \{1, \ldots , m \}$.
Then relations (4.4)' in 
$A_{aut}(\mathcal G^{\amalg n})$ become
$$\sum_{r, (r,j) \in E} X_{k(i)l(r)}
= \sum_{r, (i,r) \in E} X_{k(r)l(j)} \quad ;  \quad 
\sum_{r, (r,j) \in E} X_{l(r)k(i)}
= \sum_{r, (i,r) \in E} X_{l(j)k(r)}.$$
Then we have
$$\bullet \ \Phi_0(\sum_{r, (r,j) \in E} X_{k(i)l(r)}) = 
\sum_{r, (r,j) \in E} \nu_k(u_{ir})x_{kl} = 
\sum_{r, (i,r) \in E} \nu_k(u_{rj})x_{kl} =
\Phi_0(\sum_{r, (i,r) \in E} X_{k(r)l(j)}),$$
$$\bullet \ \Phi_0(\sum_{r, (r,j) \in E} X_{l(r)k(i)}) = 
\sum_{r, (r,j) \in E} \nu_l(u_{ri})x_{lk} = 
\sum_{r, (i,r) \in E} \nu_l(u_{jr})x_{lk} =
\Phi_0(\sum_{r, (i,r) \in E} X_{l(j)k(r)}).$$
All these computations show that $\Phi_0$ induces a $*$-homomorphism
$$\Phi : A_{aut}(\mathcal G^{\amalg n}) \longrightarrow
A_{aut}(\mathcal G) \! *_{\rm w} \! A_{aut}(X_n), 
\ \Phi(X_{k(i)l(j)}) = \nu_k(u_{ij})x_{kl}.$$
It is easily seen that $\Phi$ is a Woronowicz algebra morphism, 
that is $\Phi \circ \Delta = (\Phi \otimes \Phi) \circ \Delta$.

\medskip

We have now to construct an inverse for $\Phi$. First by Lemma 4.4
we have a $*$-homomorphism $\pi : A_{aut}(X_n) \longrightarrow 
A_{aut}(\mathcal G^{\amalg n})$ such that for $k,l \in \{1, \ldots, n \}$,
$$\pi(x_{kl}) = P^{kl} = \sum_{r=1}^n X_{k(r)l(i)}, \ 
\forall i \in \{1, \ldots, m \}.$$  
Let $\mathcal C$ be the free algebra with generators 
$(u_{ij})_{1 \leq i,j \leq m}$. For $k \in \{1, \ldots, n \}$ define an
algebra morphism 
$$\theta_0^k : \mathcal C \longrightarrow A_{aut}(\mathcal G^{\amalg n}), 
\ u_{ij} \longmapsto \sum_{l=1}^n X_{k(i)l(j)}, \ 1 \leq i,j \leq m.$$
It is immediate that $\theta_0^k$ is a $*$-homomorphism, if $\mathcal C$
is endowed with the $*$-algebra structure defined by $u_{ij}^* = u_{ij}$.
Let $i,j,j' \in \{1, \ldots m \}$. Then:
$$\bullet \ \
\theta_0^k(u_{ij}u_{ij'}) = \sum_{l,l'=1}^n X_{k(i)l(j)} X_{k(i)l'(j')}
= \sum_{l=1}^n X_{k(i)l(j)} X_{k(i)l(j')} = \delta_{jj'} \theta_0^k(u_{ij}).$$
Using Lemma 4.5, we have
$$ \bullet \ \ 
\theta_0^k(u_{ji}u_{j'i}) = \sum_{l,l'=1}^n X_{k(j)l(i)} X_{k(j')l'(i)}
= \sum_{l=1}^n X_{k(j)l(i)} X_{k(j')l(i)} = \delta_{jj'} \theta_0^k(u_{ji}).$$
$$ \bullet \ \
\theta_0^k(\sum_{j=1}^m u_{ij}) = \sum_{j=1}^m \sum_{l=1}^n X_{k(i)l(j)} = 1 = 
\theta_0^k(1).$$
Using Lemma 4.4, we have
$$\bullet \ \ \theta_0^k(\sum_{j=1}^m u_{ji}) = 
\sum_{j=1}^m \sum_{l=1}^n X_{k(j)l(i)} = \sum_{l=1}^n P_{l(i)}^{kl} 
= \sum_{l=1}^n P^{kl} = 1 = \theta_0^k(1).$$
Let $\gamma \in E$ and let $i,j \in \{1,\ldots,m \}$ with $(i,j) \not \in E$. 
Then
 \begin{align*}
& \bullet \ \ \ 
\theta_0^k(u_{s(\gamma)i}u_{t(\gamma)j}) = 
\sum_{l,l'=1}^n X_{k(s(\gamma))l(i)} X_{k(t(\gamma))l'(j)} =
\sum_{l=1}^n X_{k(s(\gamma))l(i)} X_{k(t(\gamma))l(j)} =  \\ 
& = 0 = \theta_0^k(u_{t(\gamma)j}u_{s(\gamma)i}),
\end{align*}
and using Lemma 4.5, we have
\begin{align*}
& \bullet \ \ \ 
\theta_0^k(u_{is(\gamma)}u_{jt(\gamma)}) = 
\sum_{l,l'=1}^n X_{k(i)l(s(\gamma))} X_{k(j)l'(t(\gamma))} =
\sum_{l=1}^n X_{k(i)l(s(\gamma))} X_{k(j)l(t(\gamma))} = \\ 
& = 0 = \theta_0^k(u_{jt(\gamma)}u_{is(\gamma)}).
\end{align*}
Let $\gamma, \gamma' \in E$. Using Lemma 4.5, we have
\begin{align*}
& \bullet \ \ 
\theta_0^k(u_{s(\gamma)s(\gamma')}u_{t(\gamma)t(\gamma')}) = 
\sum_{l,l'=1}^n X_{k(s(\gamma))l(s(\gamma'))} X_{k(t(\gamma))l'(t(\gamma'))} = \\
& =\sum_{l=1}^n X_{k(s(\gamma))l(s(\gamma'))} X_{k(t(\gamma))l(t(\gamma'))}  
= \sum_{l=1}^n  X_{k(t(\gamma))l(t(\gamma'))} X_{k(s(\gamma))l(s(\gamma'))}
= \theta_0^k(u_{t(\gamma)t(\gamma')}u_{s(\gamma)s(\gamma')}).
\end{align*}
Let $i,j \in \{1, \ldots, m\}$. We have:
\begin{align*}
&  \bullet \ \
\theta_0^k(\sum_{r,(r,j) \in E}u_{ir}) =
\sum_{r,(r,j) \in E} \sum_{l=1}^n X_{k(i)l(r)} =
\sum_{l=1}^n\sum_{r,(r,j) \in E} X_{k(i)l(r)} = \\
& = \sum_{l=1}^n \sum_{r,(i,r) \in E} X_{k(r)l(j)}
= \theta_0^k(\sum_{r,(i,r) \in E}u_{rj}).
\end{align*}
\begin{align*}
&  \bullet \ \
\theta_0^k(\sum_{r,(r,j) \in E}u_{ri}) =
\sum_{r,(r,j) \in E} \sum_{l=1}^n X_{k(r)l(i)} =
\sum_{l=1}^n\sum_{r,(r,j) \in E} X_{k(r)l(i)} = \\
& = \sum_{l=1}^n \sum_{r,(i,r) \in E} X_{k(j)l(r)}
= \theta_0^k(\sum_{r,(i,r) \in E}u_{jr}).
\end{align*}
Hence $\theta_0^k$ induces a $*$-homomorphism 
$\theta^k : A_{aut}(\mathcal G) \longrightarrow A_{aut}(\mathcal G^{\amalg n})$. Using the universal property of the free product, we get a $*$-homomorphism 
$\Psi_0 :  A_{aut}(\mathcal G)^{*n} \! *  \! A_{aut }(X_n)
\longrightarrow A_{aut}(\mathcal G^{\amalg n})$ such that 
$\Psi_0 \circ \nu_k = \theta^k$ and $\Psi_{0|A_{aut}(X_n)} = \pi$.
Let $k,l \in \{1, \ldots , n \}$ and $i,j \in \{1, \ldots m \}$.
Then using Lemma 4.4 and Lemma 4.5, we have:
$$\Psi_0(\nu_k(u_{ij})x_{kl}) =
\sum_{l'=1}^n \sum_{r=1}^m X_{k(i)l'(j)} X_{k(r)l(j)} =
X_{k(i)l(j)} = \Psi_0(x_{kl}\nu_k(u_{ij})).$$
In this way we get a $*$-homomorphism
$\Psi :  A_{aut}(\mathcal G) \! *_{\rm w}  \! A_{aut }(X_n)
\longrightarrow A_{aut}(\mathcal G^{\amalg n})$.
It is straightforward to check that $\Phi$ and $\Psi$ are mutually
inverse isomorphisms: this concludes the proof of Theorem 4.2.
$\square$

\bigskip

Julien Bichon

Laboratoire de Math\'ematiques Appliqu\'ees

Universit\'e de Pau et des Pays de l'Adour

IPRA, Avenue de l'universit\'e

64000 Pau, France

\smallskip

E-mail : Julien.Bichon@univ-pau.fr


\begin{thebibliography}{25}

\small{
\bibitem{[Av]} D. \textsc{Avitzour}, Free products of
$C^*$-algebras, Trans. Amer. Math. Soc. 271(2) (1982), 423-435.

\bibitem{[Ba]} T. \textsc{Banica}, Symmetries of a generic coaction,
Math. Ann. 314 (1999) 763-780.

\bibitem{[Bi]} J. \textsc{Bichon}, Quantum automorphism groups
of finite graphs, Proc. Amer. Math. Soc., to appear.

\bibitem{[Bi2]} J. \textsc{Bichon},
Galois reconstruction of finite quantum groups,
J. Algebra 230 (2000), 683-693.

\bibitem{[HR]} E. \textsc{Hewitt}, K.A. \textsc{Ross},
{\sl Abstract harmonic analysis}, Vol. 2, Springer, 1969.

\bibitem{[Il]} R. \textsc{Iltis}, Some algebraic structure
on the dual of a compact group, Canad. J. Math. 20
(1968), 1499-1510.

\bibitem{[KS]} A. \textsc{Klimyk}, K. \textsc{Schm\"{u}dgen},
{\sl Quantum groups and their representations}, 
Texts and Monographs in Physics,
Springer, 1997. 

\bibitem{[VDW]} A. \textsc{Van Daele}, S. \textsc{Wang},
Universal quantum groups, Intern. J. Math. 7 (1996), 255-263. 

\bibitem{[Voi]} D. \textsc{Voiculescu}, Symmetries of some reduced
free product $C^*$-algebras, 
in \textsl{Operator algebras and their connections with topology and ergodic theory},
Lecture Notes in Math. 1132 (1985), Springer, 556-588.

\bibitem{[W]} S. \textsc{Wang}, Free products of compact quantum groups,
Comm. Math. Phys. 167 (1995), 671-692.

\bibitem{[Wa]} S. \textsc{Wang}, Quantum symmetry groups of
finite spaces, Comm. Math. Phys. 195 (1998), 195-211.

\bibitem{[WO]} S.L. \textsc{Woronowicz}, Twisted $SU(2)$ group.
An example of noncommutative differential calculus, 
Publ. RIMS Kyoto 23 (1987), 117-181.

\bibitem{[WO1]} S.L. \textsc{Woronowicz}, Compact matrix pseudogroups,
Comm. Math. Phys. 111 (1987), 613--665.

\bibitem{[WO2]} S.L. \textsc{Woronowicz}, Tannaka-Krein duality for
compact matrix pseudogroups. Twisted $SU(N)$ groups, Invent. Math. 93 (1988),
35--76.

\bibitem{[WO2']} S.L. \textsc{Woronowicz}, A remark on compact
matrix quantum groups, Lett. Math. Phys. 21, (1991) , 35-39.

\bibitem{[WO3]} S.L. \textsc{Woronowicz}, Compact quantum groups, in 
``Sym\'etries quantiques'' (Les Houches, 1995), 
North Holland, Amsterdam, 1998, 845-884.}

\end{thebibliography}
\end{document}